\documentclass[leqno,11pt]{amsart}

\usepackage{amssymb}
\usepackage{amsmath}
\usepackage{enumerate}

\def\s{\mathbb{S}}
\def\N{\mathbb{N}}
\def\R{\mathbb{R}}

\def\Q{\mathbb{Q}}
\def\C{\mathbb{C}}

\newtheorem{theorem}{Theorem}[section]

\newtheorem{example}{Example}[section]

\newtheorem*{cora}{Corollary A}
\newtheorem*{teora}{Theorem A}

\numberwithin{equation}{section}

\begin{document}

\title[The Clifford torus as a self-shrinker]{The Clifford torus as a self-shrinker \\ for the Lagrangian mean curvature flow}

\author{Ildefonso Castro}
\address{Departamento de Matem\'{a}ticas \\
Universidad de Ja\'{e}n \\
23071 Ja\'{e}n, SPAIN} \email{icastro@ujaen.es}

\author{Ana M.~Lerma}
\address{Departamento de Matem\'{a}ticas \\
Universidad de Ja\'{e}n \\
23071 Ja\'{e}n, SPAIN} \email{alerma@ujaen.es}

\thanks{Research partially supported by a MEC-Feder grant
MTM2007-61775 and a Junta Andalucia Grant P06-FQM-01642}

\subjclass[2000]{Primary 53C42, 53B25; Secondary 53D12}

\keywords{Clifford torus, mean curvature flow, self-shrinkers,
Lagrangian surfaces.}

\date{}

\begin{abstract}
We provide several rigidity results for the Clifford torus in the
class of compact self-shrinkers for Lagrangian mean curvature
flow.
\end{abstract}

\maketitle

\section{Introduction}

An immersion $\phi : M^n \rightarrow \R^m$ of a smooth manifold
$M$ of dimension $n$ and codimension $p=m-n\geq 1$ into Euclidean
space is said to be a {\em self-shrinker} if it satisfies the
quasilinear elliptic system
\begin{equation}\label{selfshrinker}
H=-\phi^\bot
\end{equation}
where $H={\rm trace}\,\sigma$ is the mean curvature vector of the
immersion $\phi$, defined as the trace of the second fundamental
form $\sigma$, and $^\bot$ denotes the projection onto the normal
bundle of $M$. The solutions of (\ref{selfshrinker}) not only give
rise to homothetically shrinking solutions of the mean curvature
flow but also they play an interesting role in the formation of
type-1 singularities because it was proved by Huisken \cite{Hu90}
that solutions of the mean curvature flow forming such a
singularity can be homothetically rescaled so that any resulting
limiting submanifold verifies (\ref{selfshrinker}). In this way it
is expected that the understanding of the singularities of the
mean curvature flow will rely on the classification of
self-shrinkers, but this is a hard and open problem.

There are many interesting papers (for example, \cite{AL86},
\cite{Hu90}, \cite{Sm05}, \cite{CMi09}, \cite{LeSe10},
\cite{CaLi11}, \cite{DiXi11a}, \cite{DiXi11b}, \cite{ChZh11},
\cite{LiWe12} and \cite{ChPe12}) about classification and rigidity
of self-shrinkers for curves, hypersurfaces or arbitrary
codimension, under general assumptions like compactness,
completeness with polynomial volume growth, uniformly bounded
geometry, proper completeness or embeddedness. In this article we
are interested in rigidity results for {\em compact}
self-shrinkers {\em in arbitrary codimension}, emphasizing the
two-dimensional {\em Lagrangian} case. We recall that the
Lagrangian constraint is preserved by the mean curvature flow.

When $n=1$, all the solutions of (\ref{selfshrinker}) are given by
the Abresch-Langer curves \cite{AL86}. Except for the straight
line passing through the origin, their curvature is positive for
all them and the only simple closed one is the circle.

 When $m=n+1$ and $n\geq 2$,
Huisken \cite{Hu90} proved that {\em the only compact mean convex
self-shrinker is the sphere $\s^n(\sqrt n)$ of radius $\sqrt n$}.

 In higher codimension the situation becomes more complicated as
the codimension increases. A natural extension of the above
Huisken's result is the following theorem of Smoczyk \cite{Sm05}:
{\em If $M^n$ is a compact self-shrinker in $\R^m$, then $M^n$ is
spherical if and only if $|H|> 0$ and its principal normal vector
field $\nu=H/|H|$ is parallel in the normal bundle.} One easily
observes that spherical self-shrinkers coincide with minimal
submanifolds of the sphere $\s^{m-1}(\sqrt n)$ and all of them
satisfy that $|H|^2\equiv n$.

The simplest one (the totally geodesic $n$-sphere of radius $\sqrt
n $) was characterized by the following gap theorem of Cao and Li
\cite{CaLi11} (Le and Sesum \cite{LeSe10} proved first the
hypersurface case) on the squared norm of the second fundamental
form: {\em If $M^n$ is a compact self-shrinker in $\R^m$ with
$|\sigma|^2\leq 1$, then $|\sigma|^2\equiv 1 $ and $M^n$ is the
$n$-sphere $\s^{n}(\sqrt n)$ in $\R^{n+1}$.}

We observe in the above results that, in order to characterize the
sphere $\s^n(\sqrt n)$, one needs some hypothesis either on the
mean curvature or on the squared norm of the second fundamental
form. Very recently, Li and Wei in \cite{LiWe12} and  Cheng and
Peng in \cite{ChPe12} have obtained some rigidity results for
other examples of self-shrinkers under both types of assumptions:
$|H|>0$ and parallel principal normal vector field $\nu=H/|H|$ and
either lower and upper bounds of $|\sigma|^2$ or constancy of
$|\sigma|^2$.

We pay our attention to some interesting compact spherical
self-shrinkers with constant squared norm of the second
fundamental form:
\begin{example}\label{Cliford}
For any $n_1,n_2\in \N$ such that $n_1+n_2=n$, the Clifford
immersion
\[
\s^{n_1}({\sqrt n_1}) \times \s^{n_2}({\sqrt n_2}) \hookrightarrow
\R^{n+2}
\]
is a compact self-shrinker in $\R^{n+2}$ with $|\sigma|^2\equiv
2$.
\end{example}
\begin{example}\label{ProductCircles}
The product of $n$-circles
\[
\s^1 \times \stackrel{n)}{\cdots} \times \s^1 \hookrightarrow
\R^{2n}
\]
is a compact flat self-shrinker with $|\sigma|^2\equiv n$ that is
Lagrangian in $\R^{2n}\equiv\C^n$.
\end{example}
\begin{example}\label{Anciaux0}
The immersion
\[
\s^1 \times \s^{n-1} \rightarrow \C^n\equiv\R^{2n}, \quad
(e^{it},(x_1,\dots,x_n)) \mapsto \sqrt n \,
e^{it}\,(x_1,\dots,x_n)
\]
is a compact self-shrinker with $|\sigma|^2\equiv
\frac{3n-2}{n}\in [2,3)$
that is Lagrangian in $\R^{2n}\equiv\C^n$.
\end{example}

By translating the well-known results of Simon \cite{Si68}, Lawson
\cite{La69} and Chern-Do Carmo-Kobayashi \cite{CdCK78} about
intrinsic rigidity for minimal submanifolds in the unit sphere and
using some simple observations of \cite{CaLi11}, we arrive at the
following gap result for compact self-shrinkers of codimension
$p\geq 1$.

\begin{teora}\label{Th A}
Let $\phi : M^n \rightarrow \R^{n+p}$ be a compact self-shrinker
such that $|H|^2$ is constant or $|H|^2 \leq n$ or $|H|^2 \geq n$.
If
\begin{equation}\label{boundA}
|\sigma|^2\leq \frac{3p-4}{2p-3}
\end{equation}
then:
\begin{enumerate}
\item either $|\sigma|^2\equiv 1$ and $M$ is $\,\s^n(\sqrt n)$ in
$\R^{n+1}$ (i.e. $p=1$),
\item or $|\sigma|^2\equiv \frac{3p-4}{2p-3}$ and $M$ is
\begin{enumerate}
\item either $\,\s^{n_1}(\sqrt{n_1})\times\s^{n_2}(\sqrt{n_2})$, $n_1+n_2=n$, (with
$|\sigma|^2\equiv 2$) in $\R^{n+2}$ (i.e. $p=2$),
\item or the Veronese immersion (with $|\sigma|^2\equiv 5/3$) of
$\s^2(\sqrt 6)$ in $\R^5$  (i.e. $n=2$, $p=3$).
\end{enumerate}
\end{enumerate}
\end{teora}

We remark that our hypothesis on $H$ of Theorem A in the compact
case is weaker than  $|H|>0$ and parallel principal normal vector
field $\nu=H/|H|$. Our result generalizes Theorem 1.2 of
\cite{LiWe12} in the compact case.

If the dimension and the codimension of the submanifold $M$
coincide (what happens, for example, when $M$ is Lagrangian), as a
first immediate consequence of Theorem A we obtain the following
surprising characterization of the Clifford torus $\s^1 \times
\s^1 $.

\begin{cora}\label{Cor A}
Let $\phi : M^n \rightarrow \R^{2n}$ be a compact self-shrinker
with codimension $n\geq 2$ such that $|H|^2$ is constant or $|H|^2
\leq n$ or $|H|^2 \geq n$. If
\begin{equation}\label{boundAn}
|\sigma|^2\leq \frac{3n-4}{2n-3}
\end{equation}
then $n=2$, $|\sigma|^2\equiv 2$ and $M$ is the Clifford torus $\,
\s^1 \times \s^1 $ in $\, \R^4$.
\end{cora}

It is quite remarkable that the above three examples (Example 1.1,
1.2 and 1.3) coincide when $n=2$ providing precisely the Clifford
torus $\s^1 \times \s^1$ in $\R^4$. In fact, we remark that
Example 1.2 and Example 1.3 satisfy the hypothesis of Corollary A
only when $n=2$. So it can be expectable some other rigidity
results in this setting for this regular example. But the Clifford
torus $\s^1 \times \s^1$ is not isolated in the class of compact
self-shrinkers in Euclidean 4-space. In fact, it belongs to four
different families of infinitely many self-shrinkers of genus one,
that we will study in section 3:
\begin{itemize}
\item[(i)] {\em Abresch-Langer tori}, product of two
Abresch-Langer curves \cite{AL86};
\item[(ii)] {\em Anciaux tori}, defined by considering the case
$n=2 $ in Theorem 1 of \cite{An06};
\item[(iii)] {\em Lee-Wang tori}, defined by considering the case
$n=2 $ in Proposition 2.1 of \cite{LW10} and described explicitly
in Proposition 3 of \cite{CL10};
\item[(iv)] {\em Lawson tori}, described in Theorem 3 of \cite{La70}.
\end{itemize}
Thus, rigidity theorems in the family of compact (Lagrangian)
self-shrinkers in $\R^4$ are welcome and the Clifford torus is the
natural candidate for this purpose. This type of results may be
useful to try to get some progress related with the open Question
7.4 of \cite{Ne10}. Our contribution to this problem consists of
three different new characterizations of the Clifford torus,
assuming in the Lagrangian setting only one type of assumptions:
either on $H$ or on $|\sigma|^2$.

\begin{theorem}\label{Cor}
Let $\phi : M^2 \rightarrow \R^{4}$ be a compact Lagrangian
self-shrinker. If $|H|^2$ is constant or $|H|^2 \leq 2$ or $|H|^2
\geq 2$, then $M^2$ is the Clifford torus $\s^1 \times \s^1$.
\end{theorem}
As a consequence of this result we get that {\em the Clifford
torus is the only compact Lagrangian spherical self-shrinker in
$\R^{4}$}. By considering all the previous results, it seems to be
interesting the study of compact self-shrinkers with
$|\sigma|^2\leq 2$. Since there do not exist Lagrangian
self-shrinking spheres (see Theorem \ref{no spheres}), making a
subtle combination of Gauss-Bonnet Theorem with the formula
(\ref{Willmore}) that expresses the Willmore functional as a
integer multiple of the area of a compact self-shrinker, we get
the following result.

\begin{theorem}\label{Thm 4}
Let $\phi : M^2 \rightarrow \R^{4}$ be a compact orientable
Lagrangian self-shrinker. If $|\sigma|^2\leq 2$, then
$|\sigma|^2\equiv 2$ and $M$ is a torus. If, in addition, the
Gauss curvature $K$ of $M$ is non-negative or non-positive, then
$M^2$ is the Clifford torus $\s^1 \times \s^1$.
\end{theorem}


In \cite{CL10} the authors classified all Hamiltonian stationary
Lagrangian surfaces in complex Euclidean plane which are
self-similar solutions of the mean curvature flow. The Hamiltonian
stationary condition is equivalent to the vanishing of the
divergence of the tangent vector field $JH$, being $J$ the
standard complex structure of $\C^2$. Based on the above mentioned
classification, we finally deduce:

\begin{theorem}\label{Thm 1}
Let $\phi : M^2 \rightarrow \R^{4}$ be a compact self-shrinker. If
$\phi$ is a Hamiltonian stationary Lagrangian embedding, then
$M^2$ is the Clifford torus $\s^1 \times \s^1$.
\end{theorem}

After analyzing the different hypothesis in our characterizations
of the Clifford torus with the four families of self-shrinking
tori cited above, we conjecture that {\em the Clifford torus is
the only compact Lagrangian self-shrinker in $\R^4$ with
$|\sigma|^2\leq 2$}.

{\em Acknowledgements.} The first author would like to thank Knut
Smoczyk for some helpful conversations and for sharing his proof
of Theorem \ref{no spheres}.

 \vspace{0.3cm}


\section{Preliminaries}

Let $\phi : M^n \rightarrow \R^m$ be an isometric immersion of an
$n$-dimensional submanifold in Euclidean $m$-space. The mean
curvature vector $H$ of $\phi $ is given by $H={\rm
trace}\,\sigma$, where $\sigma $ denotes the second fundamental
form of $\phi $. A submanifold $M$ in $\R^m$ is called a {\em
self-shrinker} if
\begin{equation}\label{self}
H=-\phi^\bot
\end{equation}
where $^\bot$ stands for the projection onto the normal bundle of
$M$. If $^\top $ denotes projection onto the tangent bundle,  it
is easy to check that $\phi ^\top= \frac{1}{2}\nabla |\phi |^2$,
where $\nabla $ means gradient with respect to the induced metric
on $M$.

Let $\phi : M^n \rightarrow \R^m$ be a self-shrinker. Using
(\ref{self}), we get the following formula for the Laplacian of
the squared norm of $\phi$:
\begin{equation}\label{Laplacian norm}
\triangle |\phi|^2 =2 (n-|H|^2)
\end{equation}

In particular, when $M$ is compact, we obtain an interesting
relationship between the Willmore functional of $\phi $ and the
area of $M$:
\begin{equation}\label{Willmore}
\int_M |H|^2 d\mu =n \, {\rm Area}(M)
\end{equation}

\vspace{0.1cm}

 We can find many examples (see Examples 1.1, 1.2
and 1.3 in section 1) in the class of {\em spherical
self-shrinkers:}

Consider that $\phi : M^n \rightarrow \s^{m-1}(R)\subset \R^m$ is
a spherical immersion with second fundamental form  $\hat \sigma $
and mean curvature vector $\hat H$. Then $\phi $ is a
self-shrinker if and only if $\hat H =0$ and $R=\sqrt n$, that is,
$M $ is a minimal submanifold in the $(m-1)$-sphere of radius
$\sqrt n$. In this case, $H=-\phi$ and so, $|H|^2 = |\phi|^2=n$
and, in addition, it satisfies
\begin{equation}\label{spherical 2ndff}
|\sigma|^2=1+|\hat\sigma|^2
\end{equation}

\vspace{0.1cm}

We consider now a special case of codimension $n$: the Lagrangian
submanifolds, recalling that the Lagrangian constraint is
preserved by the mean curvature flow.  An immersion $\phi : M^n
\rightarrow \R^{2n}\equiv \C^n$ is said to be Lagrangian if the
restriction to $M$ of the Kaehler two-form $\omega (\,\cdot\,
,\,\cdot\,)=\langle J\cdot,\cdot\rangle$ of $\C^n$ vanishes. Here
$J$ is the complex structure on $\mathbb{C}^n$ that defines a
bundle isomorphism between the tangent and the normal bundle of
$\phi $. In particular, $\sigma (v,w)=JA_{Jv}w$, where $A$ is the
shape operator, and so the trilinear form
$C(\cdot,\cdot,\cdot)=\langle \sigma(\cdot,\cdot), J \cdot \rangle
$ is fully symmetric. On the other hand,
\begin{equation}\label{Hlagr}
H=J\nabla \beta,
\end{equation}
where $\beta:M\rightarrow \mathbb{R} /2\pi\mathbb{Z}$ is called
the Lagrangian angle map of $\phi$. In general $\beta $ is a
multivalued function; nevertheless $\alpha_H=-d\beta = \langle
JH,\cdot\rangle$ is a well defined closed 1-form on $M$ and its
cohomology class $[\alpha_H]$ is called the Maslov class of
$\phi$.

Suppose now that $\phi : M^n \rightarrow \R^{2n}\equiv \C^n$ is a
{\em Lagrangian self-shrinker}. Derivating $\phi = \phi^\top-H$ in
the direction of any tangent vector $v\in TM$ and separating
tangent and normal components, we get that
\begin{equation}\label{selfLagrT}
A_H \, v =v -\nabla_v \,\phi^\top
\end{equation}
where $A_H$ is the Weingarten endomorphism associated to $H$ and
$\nabla $ is the Levi-Civita connection of $M$, and
\begin{equation}\label{selfLagrN}
\nabla^\bot_v H = \sigma (v,\phi^\top)
\end{equation}
where $\nabla^\bot$ is the connection of the normal bundle. Using
(\ref{selfLagrN}) we obtain that
\begin{equation}\label{divJH}
{\rm div}\,JH = \langle JH, \phi^\top \rangle
\end{equation}
where div denotes the divergence operator.

We finish this section with an interesting result, first proved by
Smoczyk, about the non existence of compact orientable Lagrangian
self-shrinkers with trivial Maslov class.

\begin{theorem}\label{no spheres}
Let $\phi : M^n \rightarrow \C^n$ be a Lagrangian self-shrinker.
If $M$ is compact orientable, then $[\alpha_H]\neq 0$.
\end{theorem}
\begin{proof}
Assume that $[\alpha_H]=0$. Then there exists a globally defined
Lagrangian angle $\beta $ such that $\alpha_H=-d\beta$. Taking
into account that $\triangle \beta =-{\rm div}\,JH$, from
(\ref{divJH}) we have that $\beta$ satisfies the elliptic linear
equation  $\triangle \beta = \frac{1}{2} \langle \nabla \beta,
\nabla |\phi|^2 \rangle$. Then the maximum principle says that
$\beta $ must be constant and thus $H\equiv 0$. Since there are no
compact minimal submanifolds in Euclidean space, this is a
contradiction.
\end{proof}
As a consequence of Theorem \ref{no spheres}, {\em there do not
exist Lagrangian self-shrinkers with the topology of a sphere}.

\vspace{0.3cm}

\section{Self-shrinking tori in Euclidean 4-space}

In this section we collect the main geometric properties of four
families of self-shrinking tori all them including the Clifford
torus.

\subsection{Abresch-Langer tori}

The Abresch-Langer tori are defined simply as the product
$\Gamma_1 \times \Gamma_2$ of two closed Abresch-Langer curves
\cite{AL86}.

The curvature vector of such a curve satisfies
\begin{equation}\label{kAL}
 \overrightarrow{\kappa}_{\Gamma_i} = -\Gamma_i^\bot
\Leftrightarrow \kappa_{\Gamma_i} = \langle \Gamma_i', i \Gamma_i
\rangle, \  i=1,2
\end{equation}
where $'$ denotes derivative respect to the arclength parameter.
The only simple closed one is the unit circle and the product of
two unit circles gives obviously the Clifford torus $\s^1 \times
\s^1$. Equation (\ref{kAL}) admits a countable family of closed
noncircular solutions, which is parametrized by relatively prime
numbers $p_i$ and $q_i$ such that $p_i/q_i \in (1/2,1/\sqrt 2)$,
$i=1,2$. Moreover, if $r_i=|\Gamma_i|$, $i=1,2$, then one can
deduce from (\ref{kAL}) that $\kappa_{\Gamma_i}=\rho_i
e^{r_i^2/2}$ with $r_i^2 (1-r_i'^2)e^{-r_i^2}=\rho_i^2 $, being
$\rho_i>0$ a positive constant depending on $p_i$ and $q_i$,
$i=1,2$.

The Abresch-Langer tori are flat Lagrangian tori whose second
fundamental form and mean curvature vector satisfy
\begin{equation}\label{AandH-AL}
|\sigma|^2=|H|^2=\kappa_{\Gamma_1}^2+\kappa_{\Gamma_2}^2=\rho_1^2
e^{|\Gamma_1|^2}+\rho_2^2 e^{|\Gamma_2|^2}>0
\end{equation}
The only embedded Abresch-Langer torus is the Clifford torus.

\subsection{Anciaux tori}

Using the case $n=2$ in Theorem 1 of \cite{An06}, we can define
the Anciaux tori by the immersions $\phi_{p,q} :I \times \R
\rightarrow \C^2 $, parametrized by relatively prime numbers $p$
and $q$ such that $p/q\in (1/4,1/2)$, given by
\begin{equation}\label{An-tori}
\phi_{p,q} (t,s) = \gamma_{p,q} (t) (\cos s, \sin s)
\end{equation}
where $\gamma = \gamma_{p,q} (t)$, $t\!\in \!I \!\subset\! \R$, is
a closed curve such that its curvature satisfies the equation
\begin{equation}\label{kAn}
 \kappa_\gamma =\frac{\langle \gamma' , i\gamma
 \rangle}{|\gamma|^2}(|\gamma|^2-1)
\end{equation}
where $'$ denotes derivative respect to the arclength parameter.
It is clear that $\gamma(t)=\sqrt 2 e^{it/\sqrt 2}$ satisfies
(\ref{kAn}) providing the Clifford torus.

Following \cite{An06}
we deduce from (\ref{kAn}) that
\begin{equation}\label{kAn2}
\kappa_\gamma= \frac{E\,e^{r^2/2}(r^2-1)}{r^3}, \ r=|\gamma|, \
r^4 (1-r'^2)e^{-r^2}=E^2
\end{equation}
being $E>0$ a positive constant depending on $p$ and $q$.

The squared norm of the mean curvature vector of an Anciaux torus
is given by
\begin{equation}\label{H-An}
|H_{p,q}|^2= \frac{E^2 e^{r^2}}{r^2}
\end{equation}
and the squared norm of the second fundamental form of an Anciaux
torus is given by
\begin{equation}\label{A-An}
|\sigma_{p,q}|^2=\frac{E^2 e^{r^2}}{r^6}(r^4-2r^2+4)
\end{equation}
Every Anciaux torus is Lagrangian but the only embedded one is the
Clifford torus by Theorem 3 in \cite{An06}.

\subsection{Lee-Wang tori}

From \cite{CL10} and \cite{LW10}, we define the Lee-Wang tori
$\mathcal{T}_{m,n}$ by the doubly-periodic immersions
$\Psi_{m,n}:\mathbb{R}^2 \rightarrow \mathbb{C}^2$, $m,n\in\N, \,
(m,n)=1, \, m\leq n $, given by
\begin{equation}\label{LW-tori}
\Psi_{m,n}(s,t)=\sqrt{m+n} \left( \frac{1}{\sqrt{n}}\,\cos s\,e^{i
\sqrt{\frac{n}{m}}t} , \frac{1}{\sqrt{m}}\,\sin s\,e^{i
\sqrt{\frac{m}{n}}t} \right)
\end{equation}
The Clifford torus corresponds to $\mathcal{T}_{1,1}$, since
$\Psi_{1,1}(s,t)=\sqrt{2} e^{i t} ( \cos s  , \sin s\, )$.

The Lee-Wang tori are Hamiltonian stationary Lagrangian tori
satisfying
\begin{equation}\label{normLW}
\frac{m+n}{n}\leq |\Psi_{m,n}|^2=\frac{m+n}{mn}(m\cos^2 s +n\sin^2
s) \leq \frac{m+n}{m}
\end{equation}
After a straightforward computation, the squared norm of the mean
curvature vector of a Lee-Wang torus is given by
\begin{equation}\label{H-LW}
1 < \frac{m+n}{n} \leq |H_{m,n}|^2=\frac{m+n}{n\cos^2 s + m \sin^2
s}\leq  \frac{m+n}{m}
\end{equation}
and the squared norm of the second fundamental form and the Gauss
curvature of a Lee-Wang torus verify
\begin{equation}\label{A-LW}
\frac{3m^2+n^2}{n(m+n)} \leq |\sigma_{m,n}|^2 \leq
\frac{m^2+3n^2}{m(m+n)}
\end{equation}
and
\begin{equation}\label{K-LW}
-\frac{n(n-m)}{m(m+n)} \leq K_{m,n} \leq \frac{m(n-m)}{n(m+n)}
\end{equation}

The only embedded Lee-Wang torus is the Clifford torus by
Proposition 3 in \cite{CL10}.

\subsection{Lawson tori}
We define (cf.\ \cite{La70}) the Lawson tori $\mathcal{T}_\alpha$
by the doubly-periodic immersions $\Phi_\alpha: \mathbb{R}^2
\rightarrow \mathbb{C}^2$, $\alpha \in \Q, \, \alpha \geq 1 $,
given by
\begin{equation}\label{Law-tori}
\Phi_\alpha(x,y)=\sqrt 2 \left( \cos x \, e^{i\alpha y} , \sin x
\, e^{iy}  \right)
\end{equation}
The Clifford torus corresponds to $\mathcal{T}_1$, since
$\Phi_1(x,y)=\sqrt 2 e^{iy} ( \cos x   , \sin x )$.

The Lawson tori are spherical tori whose squared norm of the
second fundamental form is given by
\begin{equation}\label{A-Law}
 1+\frac{1}{\alpha^2} \leq |\sigma_\alpha|^2 = 1 + \frac{\alpha^2}{(\alpha^2 \cos^2 x + \sin^2 x)^2}
 \leq 1+\alpha^2
\end{equation}
The Gauss curvature of a Lawson torus verifies
\begin{equation}\label{K-Law}
1-\alpha^2 \leq K_\alpha \leq 1-\frac{1}{\alpha^2}
\end{equation}

The only embedded Lawson torus is the Clifford torus, which is
also the only Lagrangian in this family.

\vspace{0.3cm}


\section{Proof of the results}
In this section we prove the results stated in section 1.

{\em Proof of Theorem A.}

We first integrate (\ref{Laplacian norm}) using that $M$ is
compact, obtaining $ 0=\int_M (n-|H|^2)d\mu$. By the hypothesis on
the mean curvature $H$ of $\phi $ we conclude $|H|^2\equiv n$.
Using (\ref{Laplacian norm}) again, we deduce that $|\phi|^2$ is
harmonic, so it must be constant. Using section 2, we have that
$\phi $ is a minimal submanifold in the sphere $\s^{n+p-1}(\sqrt n
)$.

In addition, from (\ref{spherical 2ndff}) and (\ref{boundA}), we
obtain that its second fundamental form $\hat \sigma$ satisfies
\begin{equation}\label{boundAhat}
|\hat \sigma|^2\leq \frac{p-1}{2p-3}.
\end{equation}
We now recall the well known results of \cite{Si68}, \cite{La69}
and \cite{CdCK78}, about intrinsic rigidity for minimal
submanifolds in the unit sphere, which can be summarized in the
folllowing way:

{\em If $M^n$ is a compact minimal submanifold of $\,\s^{n+q}$
with second fundamental form $\hat \sigma$ such that $|\hat
\sigma|^2 \leq n/(2-1/q)$, then either $|\hat \sigma|^2 \equiv 0$
and $M^n$ is $\,\s^n$, or $|\hat \sigma|^2 \equiv n/(2-1/q)$ and
$M^n$ is $\,\s^k(\sqrt{\frac{k}{n}})\times
\s^{n-k}(\sqrt{\frac{n-k}{n}})$ in $\,\s^{n+1}$, $1\leq k \leq
n-1$, or the Veronese immersion of $\,\s^2(\sqrt 3)$ in $\,\s^4$.}

Up to a dilation in $\R^{n+q+1} \supset \s^{n+q} $ of ratio $\sqrt
n$, we can rewrite it as follows:

{\em If $M^n$ is a compact minimal submanifold of
$\,\s^{n+q}(\sqrt n )$ with second fundamental form $\hat \sigma$
such that $|\hat \sigma|^2 \leq q/(2q-1)$, then either $|\hat
\sigma|^2 \equiv 0$ and $M^n$ is $\,\s^n (\sqrt n)$, or $|\hat
\sigma|^2 \equiv q/(2q-1)$ and $M^n$ is $\,\s^k(\sqrt k)\times
\s^{n-k}(\sqrt{n-k})$ in $\,\s^{n+1}(\sqrt n)$, $1\leq k \leq
n-1$, or the Veronese immersion of $\,\s^2(\sqrt 6)$ in
$\,\s^4(\sqrt 2)$.}

Taking $q=p-1$, thanks to (\ref{boundAhat}), we make use of the
above result to finish the proof of Theorem A.

\vspace{0.2cm}

{\em Proof of Theorem~\ref{Cor}.}

The same argument of the first part of the proof of Theorem A
implies that $\phi$ is spherical. So we know from section 2 that
$|\phi|^2=|H|^2\equiv 2$ and Theorem 1.1 in \cite{Sm05} says that
$\nabla^\bot H =0$. Since $\phi$ is Lagrangian we have that $JH$
is a non-null parallel tangent vector field on $M$. We use then
Theorem 3 in \cite{Ur87} to deduce that $M$ must be necessarily a
standard torus, product of two circles. Using finally that $\phi $
is a self-shrinker, $M$ only can be the Clifford torus.

\vspace{0.2cm}

{\em Proof of Theorem~\ref{Thm 4}.}

We use Gauss-Bonnet Theorem in the Gauss equation of $\phi$
\begin{equation}\label{Gauss eq}
2K=|H|^2-|\sigma|^2
\end{equation}
obtaining
\begin{equation}\label{formula}
8\pi (1-{\rm gen}(M)) =2\int_M K\,d\mu=\int_M
(|H|^2-|\sigma|^2)d\mu = \int_M (2-|\sigma|^2)d\mu,
\end{equation}
the last equality thanks to (\ref{Willmore}).

Theorem \ref{no spheres} implies that $M$ can not be a sphere.
Hence the hypothesis $|\sigma|^2 \leq 2$ in (\ref{formula}) says
that $M$ is a torus with $|\sigma|^2\equiv 2$. If, in addition,
$K\geq 0$ or $K\leq 0$, we deduce from (\ref{Gauss eq}) that
$|H|^2\geq 2$ or $|H|^2\leq 2$. Then Theorem \ref{Cor} gives that
$M$ is the Clifford torus.

\vspace{0.2cm}

{\em Proof of Theorem~\ref{Thm 1}.}

In Corollary 1 of \cite{CL10}, the authors proved that the
Lee-Wang tori ${\mathcal T}_{m,n}$ are the only compact orientable
Hamiltonian stationary Lagrangian self-shrinkers. Since the
Clifford torus ${\mathcal T}_{1,1}$ is the only embedded in this
family, we finish the proof using that a Klein bottle does not
admit a Lagrangian embedding in $\C^2$ (see \cite{Nm09}).

\vspace{0.3cm}


\end{document}